\documentclass[12pt]{amsart}
\usepackage{amssymb}
\ifx\mathrm\undefined\let\mathrm\rm\fi
\ifx\mathbf\undefined\let\mathbf\bf\fi
\ifx\mathfrak\undefined\let\mathfrak\frak\fi
\ifx\mathcal\undefined\let\mathcal\cal\fi
\ifx\mathbb\undefined\let\mathbb\Bbb\fi
\ifx\emph\undefined\let\emph\it\fi

\newcommand{\Z}{{\mathbb Z}}

\newcommand{\R}{{\mathbb R}}

\newcommand{\Ref}[1]{{(\ref{#1})}}
\newcommand{\be}{\begin{displaymath}}
\newcommand{\ee}{\end{displaymath}}
\newcommand{\bea}{\begin{eqnarray*}}
\newcommand{\eea}{\end{eqnarray*}}

 \let\eps\varepsilon \let\epsilon\eps

\let\al\alpha
\let\bt\beta
\let\gm\gamma
\let\Gm\Gamma

\let\ka\kappa
\let\la\lambda 
 
 \let\phi\varphi

\let\si\sigma  
\let\tht\theta  \let\thi\vartheta

\newcommand{\half}{\frac12}
\newcommand{\bean}{\begin{eqnarray}}
\newcommand{\eean}{\end{eqnarray}}

\newcommand{\vs}{\vspace{.5\baselineskip}}
\newtheorem%
{thm}{Theorem}%[section]
\newtheorem%
{proposition}[thm]{Proposition}
\newtheorem%
{lemma}[thm]{Lemma}
\newtheorem%
{lemmadef}[thm]{Lemma-Definition}
\newtheorem%
{corollary}[thm]{Corollary}
\newtheorem%
{conjecture}[thm]{Conjecture}
\title[Elliptic Selberg Integrals]
{Elliptic Selberg Integrals}
\author[G. Felder, L. Stevens, and A. Varchenko]
{G. Felder $^{\,\star}$,
L. Stevens$^{\, \diamond}$, and
A. Varchenko$^{\,\diamond}$}
\begin{document}
\maketitle

\begin{center}
{\it
$^\star$ Departement Mathematik, ETH-Zentrum, 8092 Z\"urich, Switzerland,

felder@math.ethz.ch

\medskip

$^\diamond$ Department of Mathematics, University of North Carolina
at Chapel Hill,

Chapel Hill, NC 27599-3250, USA,

stevens@math.unc.edu, av@math.unc.edu}
\end{center}

%\vsk1.5>
\centerline{March, 2001}
%\vsk1.8>

%\begin{abstract}

%\end{abstract}

%\vsk>
%\vsk0>
The Selberg integral is
\bean
B_{p}(\al,\bt,\gm)=
\int_{\Delta_p}\prod_{j=1}^{p}t_{j}\,^{\al-1}(1-t_{j})^{\bt-1}
\prod_{0\leq j<k\leq 1}(t_{j}-t_{k})^{2\gm},
\notag
\eean
where $\Delta_{p}=\{t\in\R^{p}\,|\,0\leq t_{p}\leq\cdots\leq t_{1}\leq 1\}.$
The Selberg integral is a generalization of the beta function.  It can
be calculated explicitly,
\bean
B_{p}(\al,\bt,\gm)=
{1\over p!}
\prod_{j=0}^{p-1}
\frac{\Gm(1+\gm+j\gm)
\Gm(\al+j\gm)
\Gm(\bt+j\gm)}
{\Gm(1+\gm)
\Gm(\al+\bt+(p+j-1)\gm)}.
\notag
\eean
The Selberg integral has many applications, see \cite{A1,A2,As,D,DF1,DF2,M,S}.  In this paper, we
present an elliptic version of the Selberg integral.

Let $\thi_{1}(t,\tau)$ be the first Jacobi theta function \cite{WW},
\bean
\thi_{1}(t,\tau)=-\sum_{j\in\Z}
e^{\pi i(j+\half)^2\tau+2\pi i(j+\half)(t+\half)}.
\notag
\eean

Introduce special functions 
\bean
\si_{\la}(t,\tau)=
\frac{\thi_{1}(\la-t,\tau)\thi_{1}'(0,\tau)}
{\thi_{1}(\la,\tau)\thi_{1}(t,\tau)},
\qquad
\rho(t,\tau)=\frac{\thi_{1}'(t,\tau)}{\thi_{1}(t,\tau)},
\qquad
E(t,\tau)=\frac{\thi_{1}(t,\tau)}{\thi_{1}'(0,\tau)}.
\notag
\eean
Here $'$ denotes the derivative with respect to the first argument.
%The functions $E(t,\tau)^a$ and $\si_{\la}(t,\tau)$ are
%elliptic analogs of the functions $t^a$ and $1/t$, respectively.

Let $\ka\geq2$ be an integer.  The theta functions
\bean
\tht_{\ka,m}(\la,\tau)=
\sum_{j\in\Z}e^{2\pi i\ka(j+\frac{m}{2\ka})^{2}\tau+2\pi
  i\ka(j+\frac{m}{2\ka})\la},\qquad m\in\Z/2\ka\Z,
\notag
\eean
form a basis of the theta functions of level $\ka$. 

For a positive integer $p$, the elliptic Selberg integral
$I_{p}(\la,\tau)$ is the integral, 
\bean
I_{p}(\la,\tau)=J_{p}(\la,\tau)+(-1)^{p+1}J_{p}(-\la,\tau),
\notag
\eean
where
\bean
J_{p}(\la,\tau)
=
\int_{\Delta_{p}} 
\prod_{j=1}^p E(t_{j},\tau)^{-\frac{p}{p+1}}
\prod_{1\leq j<k\leq
  p}E(t_{j}-t_{k},\tau)^{\frac{1}{p+1}}\,&\times\notag
\\
\prod_{j=1}^p\si_{\la}(t_{j},\tau)
\tht_{2(p+1),p+1}\left(\la+\frac{1}{p+1}\sum_{j=1}^p
t_{j},\tau\right)&dt_{1}\dots dt_{p}\,\,.  
\notag
\eean
The branch of the logarithm is chosen in such a way that
arg $(E(t,\tau))\to 0$ as $t\to 0^{+}$, and the integral is understood as a
natural analytic continuation.\footnote{To define the analytic
  continuation, we replace the exponential $-{p\over p+1}$
  by $a$, and consider the analytic continuation with respect to $a$
  from the region where $a$ is positive.} 

\begin{thm}
We have
\bean\label{TM}
 I_{p}(\la,\tau)=c_{p}\,B_p\left(\half+{1\over2(p+1)},-{p\over p+1}, {1\over
   2(p+1)}\right)\,
\thi_{1}(\la,\tau)^{p+1}
\eean
where
\bean
c_{p}=-(2\pi)^{{p\over2}}\,e^{\pi i{p\over p+1}}
\,e^{-\pi i{p+2\over4}}\,\prod_{j=1}^{p}\left(1-e^{-\pi i{j\over p+1}}\right).  
\notag
\eean
\end{thm}
The theorem is a generalization of theorem $13$ in \cite{FV1}.  The
proof is based on the following remarks.  Consider the heat equation
\bean
4\pi i(p+1)\frac{\partial u}{\partial\tau}(\la,\tau)=
\frac{\partial^{2}u}{\partial\la^{2}}(\la,\tau)+
p(p+1)\rho'(\la,\tau)u(\la,\tau).
\notag
\eean
It is known that this equation has a one dimensional space of
solutions $u(\la,\tau)$ which are holomorphic theta functions of level
$2(p+1)$,
\bean
u(\la+2,\tau)=u(\la,\tau), \qquad
u(\la+2\tau,\tau)=e^{-  4\pi i(p+1)(\la+\tau)}u(\la,\tau)
\notag
\eean
and Weyl skew-symmetric, $u(-\la,\tau)=(-1)^{p+1}u(\la,\tau)$, see
\cite{FV1,FV2}.  The space is called the space of conformal blocks.
Clearly the right hand side of \Ref{TM} has these properties.  According to
\cite{FV1}, the left hand side of \Ref{TM} also has these properties.  Thus the
two functions are proportional.  The coefficient of proportionality is
easily calculated in the limit $\tau\to i\infty$.

\end{document}